# Hirotaka's problem 028

De Boeck I.

## Introduction

In (Majewski, 2020), the authors discuss some sangaku's made by the Japanese mathematician Hirotaka Ebisui. Most of them only require a basic knowledge of plane Euclidean geometry to be solved. One of these is called problem HI 028. In what follows, we will give a proof and make some extra observations.

## The problem

Given two perpendicular lines and two circles tangent to both lines. The circles are on the same side of one of the lines and on different sides of the other line.

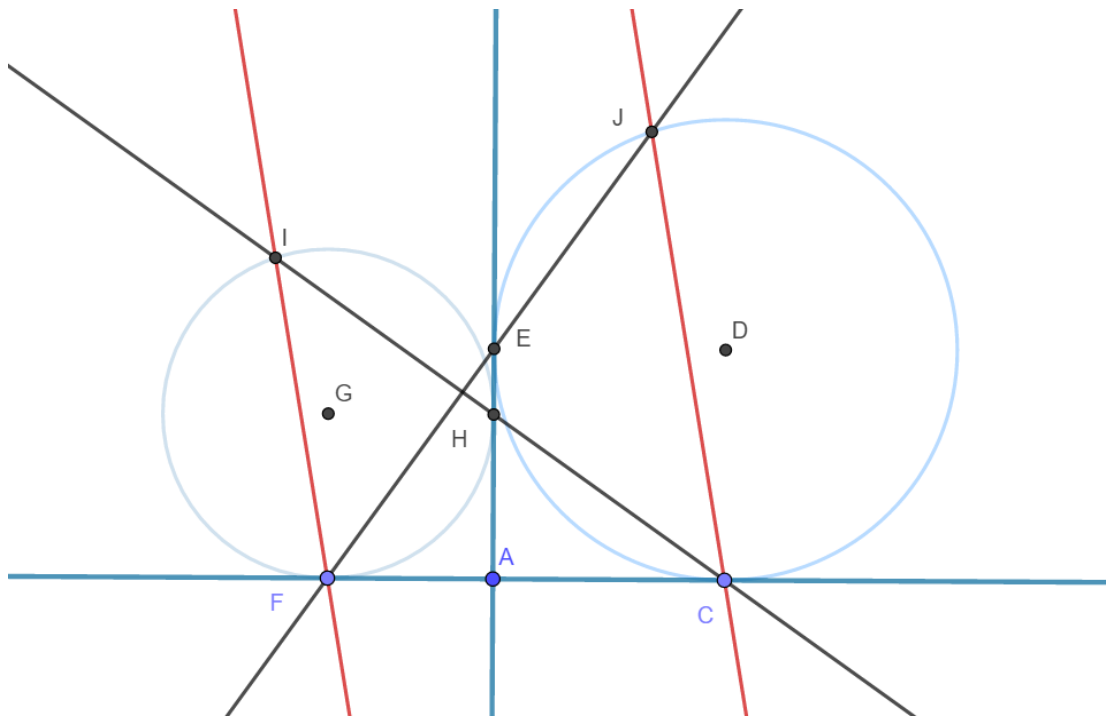

*Figure 1*

Here we have the circles c and c' with midpoints G and D respectively and the perpendicular lines AE and AC. The circle c' is tangent to AC in C and to AE in E. The circle c touches the line AC in F and the line AE in H. Note that HE is one of the internal common tangents of the circles c and c' and FC one of the external common tangents. So one of the external common tangents is perpendicular to one of the internal common tangents.

We connect one of the tangent points on the internal common tangent of one of the circles with the tangent point of the external common tangent on the other circle, e.g. E with F and H with C. These lines intersect the circles in two other points: I and J. The objective is to prove that the lines IF and JC are parallel.

But there is more in the picture. We can also prove that the lines IC and FJ are perpendicular and that IJ is the other common external tangent of the circles. Moreover, if we indicate the other two intersection points of the lines IC and FJ with the circles, we find the other internal common tangent, LM in the next figure. This internal common tangent is also perpendicular to the external common tangent IJ.

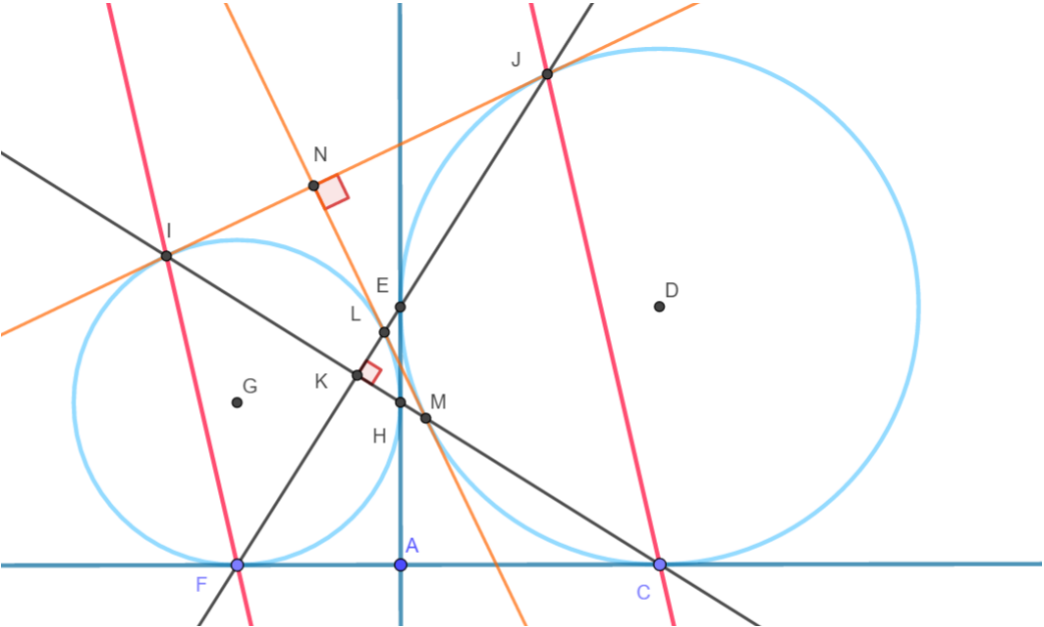

*Figure 2*

We will try to prove these statements.

## Proof of the statements

Let K be the intersection point of the lines HC and FE. First note that the angles $K\hat{J}C$ and $C\hat{I}F$ are equal to 45° since they are circumferential angles on the same arc as the 90° center angles $E\hat{D}C$ and $H\hat{G}F$.

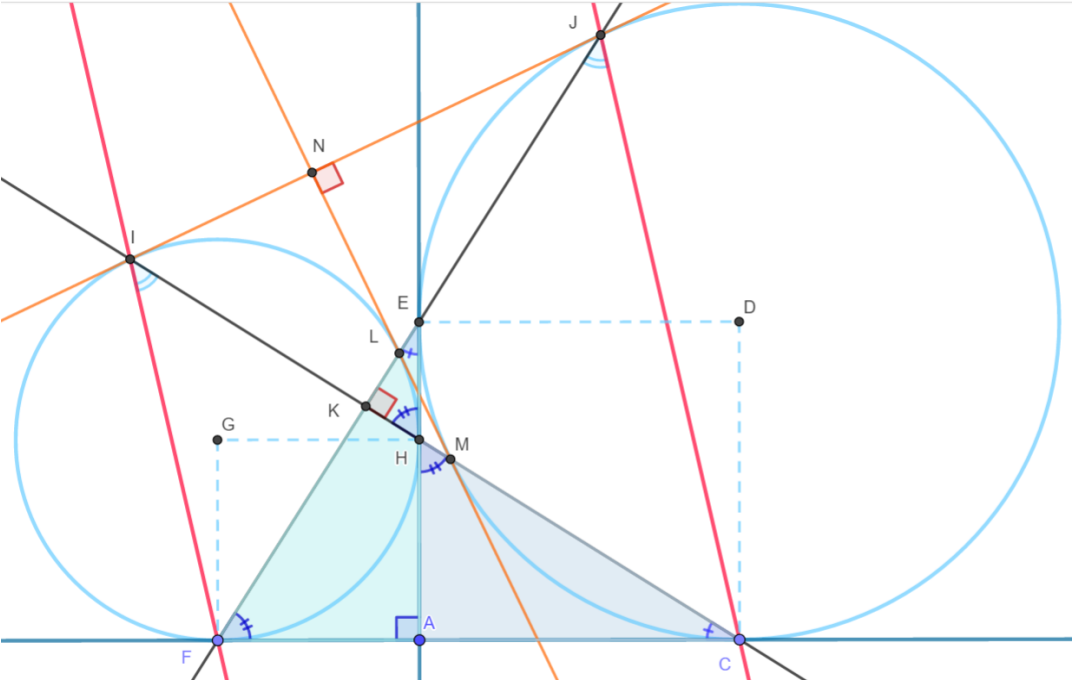

*Figure 3*

Since the triangles HCA and FEA are congruent (|AC| = |AE|, $C\hat{A}H = E\hat{A}F = 90°$, |AH| = |AF|), $F\hat{E}A = H\hat{C}A = 90° - A\hat{H}C = 90° - E\hat{H}K.$ Hence, in triangle EKH the angle $E\hat{K}H = 90°$.

Now this implies that the angle $K\hat{C}J$ in triangle JKC is 45°, as is $I\hat{F}K$ in triangle IKF. It then follows that $F\hat{I}K = K\hat{C}J$, so JC // IF.

This proves the statement made by Hirotaka and the statement that the lines IC and JF are perpendicular.

In order to see that IJ is the other common external tangent, we note that both triangles JKC and IKF are isosceles and right angled. The triangles KCF and KJI are then congruent since |KC| = |KJ|, the angles in K are opposite angles and |KF| = |KI|.

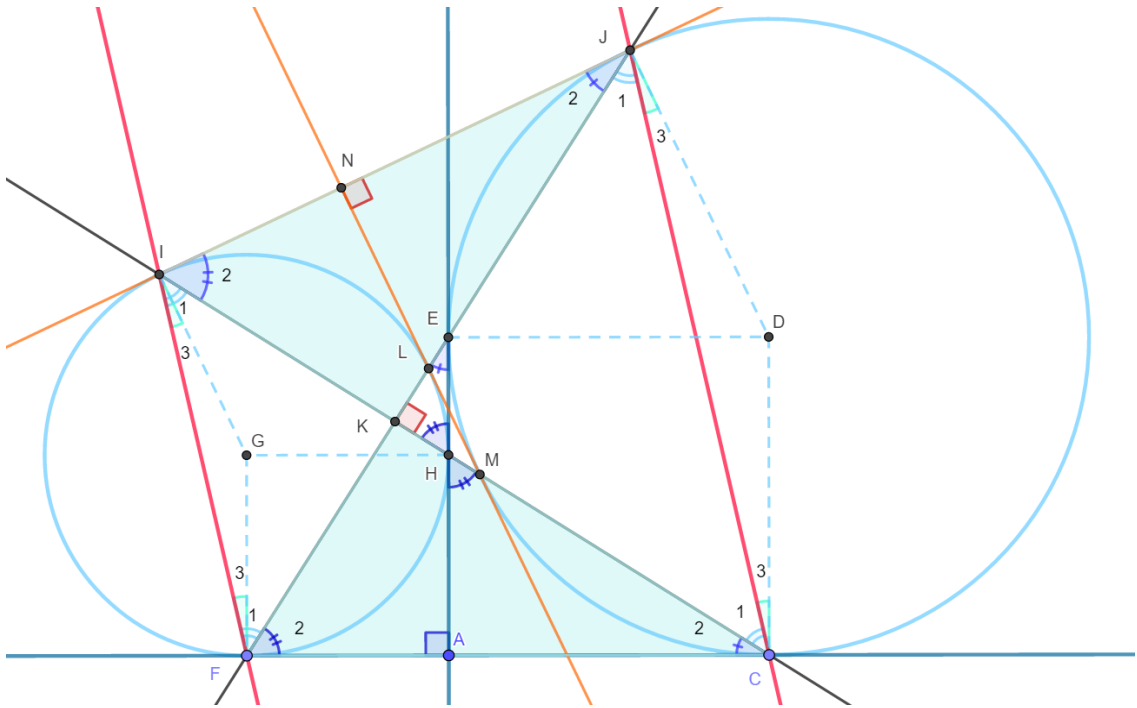

*Figure 4*

The angle between the lines IJ and JD is 90° because each of its parts is equal to the angle in C with the corresponding number (see figure 4): those with number one are 45°, the angles with number two are corresponding angles in the congruent triangles JKI and CKF and the angles numbered by 3 are the basic angles of the isosceles triangle JDC. Moreover, the sum of the angles in C equals 90° because FC is tangent in C to the circle c'.

A similar reasoning can be made in the points I and F, but now we have to consider the sum of the angles numbered by 1 and 2 decreased with the angles numbered by 3. So we proved the second claim.

For the third claim, we first note that $\hat{M}_2 = \hat{M}_3 = 90° - K\hat{L}M = 90° - N\hat{L}J = \hat{J}_2 = \hat{C}_2$ (see figure 5). Moreover, since the triangle MDC is isosceles, $\hat{M}_1 = \hat{C}_1 + \hat{C}_3$, which means that
$\hat{M}_1 + \hat{M}_2 = \hat{C}_1 + \hat{C}_3 + \hat{C}_2 = 90°$. Hence LM is perpendicular to MD.

*Figure 5*

In a similar way, $\hat{L}_1 = \hat{F}_1 - \hat{F}_3$, which means that
$\hat{L}_1 + \hat{L}_2 = \hat{F}_1 - \hat{F}_3 + \hat{F}_2 = 90°$. Hence LM is perpendicular to GL. It follows that LM is the other common internal tangent of the two circles and that $\hat{L}_1 = \hat{J}_2$.

*Figure 6*

### First conclusion

As a conclusion, we can make the following statements. If two mutually external circles with no points in common have a perpendicular external and internal common tangent, then the other internal and external common tangent are also perpendicular. Moreover, the tangent points lie on two perpendicular lines and the lines connecting two outer tangent points of the same circle are parallel.

Note that this last conclusion also follows from the fact that lines connecting the tangent points of the external common tangents of two mutually external circles with no common points are always parallel. Indeed, if we consider the homothety with center the intersection point X of the two tangents and coefficient $\frac{|XI|}{|XJ|}$, the image of the line IF is the line JC. Hence, these two lines are parallel.

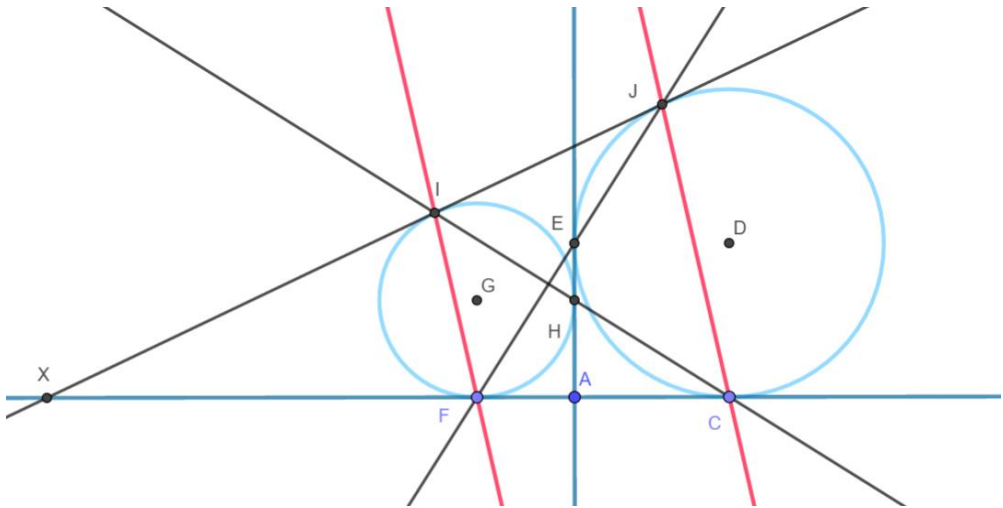

*Figure 7*

### Reverse statement

We can also investigate the reverse statement: if all the tangent points of two mutually external circles with no points in common are on two lines, these lines are perpendicular and the internal and external common tangents are two by two perpendicular.

We start with two mutually external circles c and c' with no common points. The lines IJ and FC are the external common tangents and the lines LM and EH the internal common tangents. We suppose that the tangent points J, E, L and F are collinear as well as the tangent points I, H, M and C.

### Symmetry

First, we note that the line connecting the centers of the two circles is an axis of symmetry. The reflection of one internal (external) common tangent to this axis yields the other internal (external) common tangent and vice versa. So the reflection of one of the lines connecting tangent points yields the other line connecting tangent points. Hence, the intersection point K of these two lines lies on the line connecting the centers of the circles. We will use this symmetry in the argument.

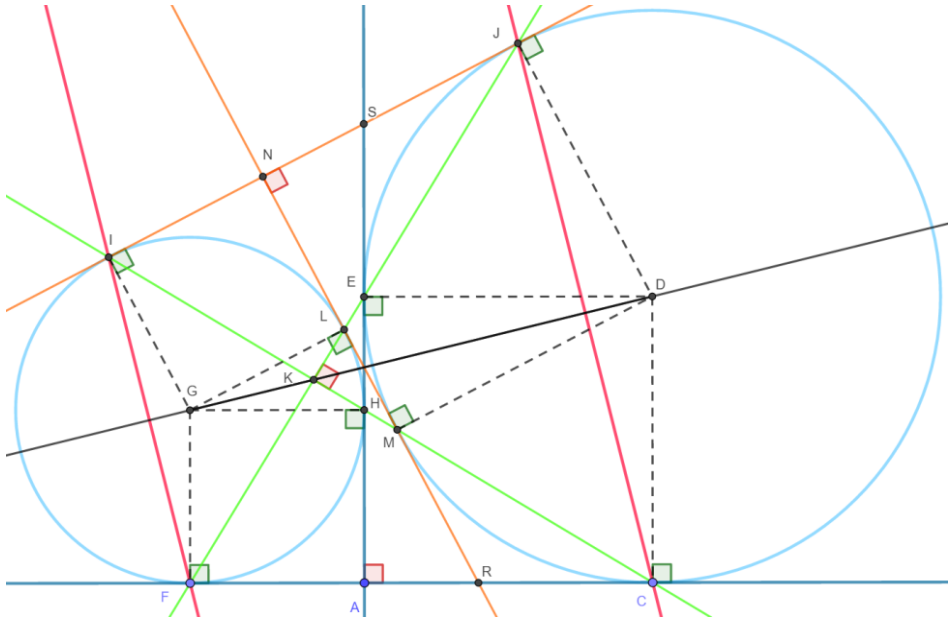

*Figure 8*

**Proof that the internal and external common tangents are perpendicular**

First note that |AE| = |AC| (property of the tangent line segments or consider the two congruent triangles ADC and ADE) and |AF| = |AH| (congruent triangles GHA and GFA).

Second, we look at the angles. Because the triangle MDC is isosceles, the angles $D\widehat{M}C$ and $D\widehat{C}M$ are equal. Hence, the complementary angles $R\widehat{M}C$ and $M\widehat{C}R$ are also equal. By symmetry, it follows that $F\widehat{E}A = J\widehat{E}S = R\widehat{M}C = M\widehat{C}R$.

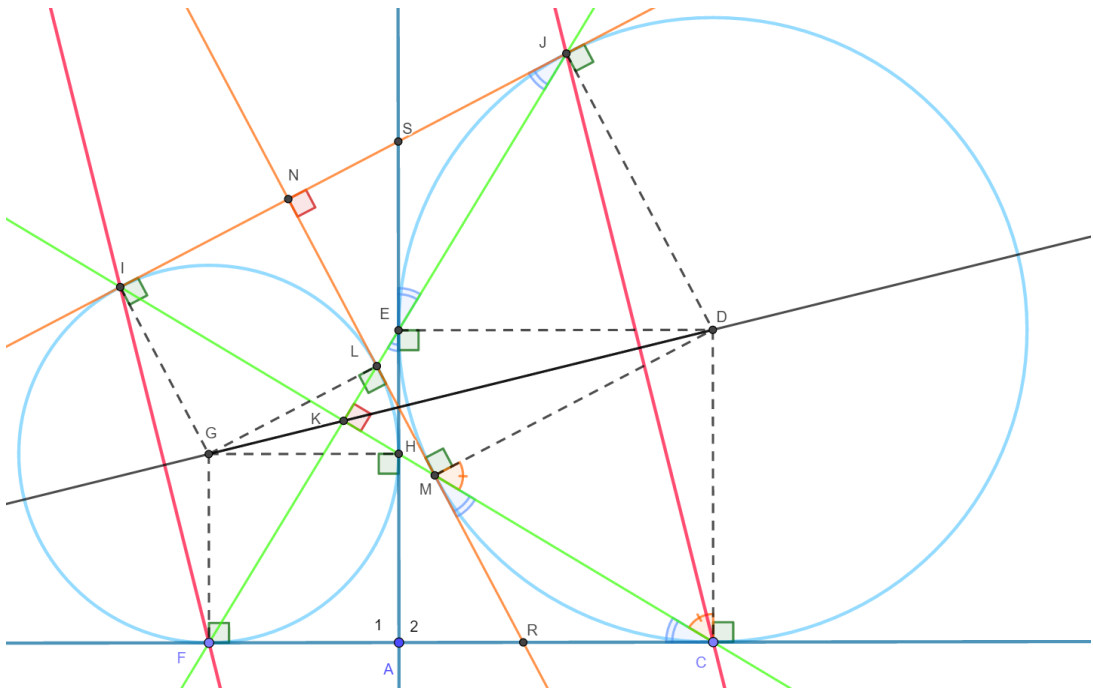

*Figure 9*

We combine these findings to state that the triangles EAF and CAH are congruent. By consequence, the angles in A are equal and thus both 90°.

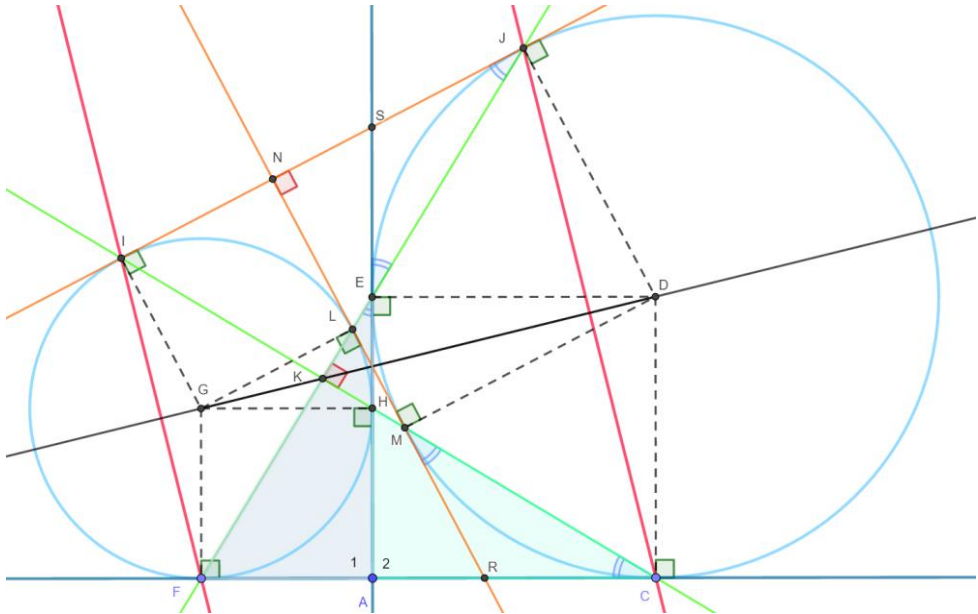

*Figure 10*

Symmetry then yields that the angles in N are also 90°.

Moreover, the quadrilateral EDCA is a square which implies that the angle $K\hat{J}C$ is a circumferential angle on a quarter circle. Thus $K\hat{J}C = 45°.$ By symmetry $J\hat{C}K = 45°$, which means that $J\hat{K}C = 90°$.

## Conclusion

Problem HI 028 of Hirotaka led us to the following theorem: two mutually external circles with no common points have al tangent points lying on two lines if and only if the internal and external common tangents are two by two perpendicular. Moreover, the lines on which the tangent points lie, are perpendicular.